\documentclass[12pt,a4paper]{article}
\usepackage[latin2]{inputenc}
\usepackage{amsmath}
\usepackage{amsfonts}
\usepackage{amssymb}
\usepackage{graphicx}
\usepackage[left=2.00cm, right=2.00cm, top=2.50cm, bottom=2.50cm]{geometry}
\usepackage[T1]{fontenc}

\def\Q{{\mathbb Q}}
\def\Z{{\mathbb Z}}

\newtheorem{lemma}{Lemma}
\newtheorem{theorem}[lemma]{Theorem}

\title{
Calculating relative power integral bases \\
in totally complex quartic extensions of totally real fields
}
\author{
Istv\'{a}n Ga\'{a}l,\; \\
{\small University of Debrecen, Mathematical Institute} \\
{\small H--4002 Debrecen Pf.400., Hungary,} 
{\small e--mail: gaal.istvan@unideb.hu},
}

\begin{document}
\baselineskip=17pt

\maketitle
\thispagestyle{empty}

\renewcommand{\thefootnote}{\arabic{footnote}}
\setcounter{footnote}{0}

\noindent
Mathematics Subject Classification: Primary 11Y50; 
Secondary 11R04, 11D57, 11D59.\\
Key words and phrases: power integral basis; calculating solutions of
index form equations;
relative quartic extensions;
unit equation; Thue equation

\begin{abstract}
Some time ago we extended our monogenity investigations and
calculations of generators of power integral bases to the relative case, cf.
\cite{book}, \cite{gsz13}, \cite{grsz}. Up to now we 
considered (usually totally real) extensions of complex quartic fields.

In the present paper we consider power integral bases in
relative extensions of totally real fields. 
Totally complex quartic extensions of totally real number fields
seems to be the most simple case, that we detail here.
As we shall see, even in this case we have to overcome several
unexpected difficulties, which we can, however solve by properly 
(but not trivially) adjusting standard methods.
We demonstrate our general algorithm on an explicit example.
We describe how the general methods for solving relative index form
equations in quartic relative extensions are modified in this case.
As a byproduct we show that relative 
Thue equations in totally complex extensions of totally real fields
can only have small solutions, and we construct a special 
method for the enumeration of small solutions of special unit equations.
These statements can be applied to other diophantine problems,
as well.
\end{abstract}

\newpage

\section{Introduction}

Monogenity of number fields have an extensive literature cf. \cite{nark}, \cite{book}.
A number field $K$ of degree $n$ is monogenic if its ring of integers $\Z_K$ is a 
simple ring extension of $\Z$, that is if there exists $\alpha\in\Z_K$
such that $\Z_K=\Z[\alpha]$. In this case $(1,\alpha,\ldots,\alpha^{n-1})$
is an integral basis of $K$, called power integral basis.

In the relative case, if $K$ is an extension field of $L$ with $[K:L]=n$,
$K$ is (relatively) monogenic over $L$ if $\Z_K$ is a simple ring extension
of the ring of integers $\Z_L$ of $L$, that is if there exists 
$\alpha\in\Z_K$ with $\Z_K=\Z_L[\alpha]$ (see \cite{book}).
In this case $(1,\alpha,\ldots,\alpha^{n-1})$ is a relative integral
basis of $K$ over $L$ called relative power integral bases. 
If $K$ is monogenic, then it is relatively monogenic over its subfield $L$,
cf. \cite{grsz}.

There are efficient algorithms for calculating generators of
power integral bases in lower degree number fields, up to degree 5
(see \cite{book}), but there is no general efficient algorithm
for number fields of arbitrary degrees. These algorithms are especially
efficient for cubic and quartic number fields when the problem can be
reduced to the resolution of (one or more) Thue equations.

Some time ago we extended our investigations to the relative case
cf. \cite{book}, \cite{gsz13}, \cite{grsz}. At first we studied
extensions of complex quadratic fields. In the present paper
we consider relative extensions of totally real number fields
(of arbitrary degrees). The most simple case seems to be 
the relative quartic extensions, when we can follow the arguments of 
\cite{gprel4}, reducing the calculation of generators of 
relative power integral bases to solving relative Thue equations. 
Especially, if $K$ is a totally complex quartic extension of the
totally real number field $L$, then, at first glance, the calculation
seems to be easy. Trying to perform the calculation for an explicite
example we meet some unexpected problems. The resolution of
the cubic Thue equation (according to \cite{gprel4}) is not trivial,
even the standard enumeration algorithm has to be adjusted properly.
Also, we observe that if $K$ is a totally complex quartic extension of
a totally real number field $L$, then the quartic relative Thue equation
to be solved by \cite{gprel4} is trivial, having only small solutions
(cf. Theorem \ref{thth}). 

These observations might well be applied also in several other calculations.

We  illustrate our paper by an explicit example that gives a direct insight
into the calculations.

\section{The general scheme for relative quartic extensions}

We recall the general method of \cite{gprel4} which reduces
calculation of generators of relative power integral bases 
in quartic relative extensions to the resolution
of relative Thue equations.

In the following we denote by $\Z_K$ the ring of algebraic integers of any
number field $K$.

Let $M$ be a totally real number field of degree $m$. 
Assume $M=\Q(\mu)$ with $\mu\in\Z_M$.
Let $K=M(\xi)$ be a totally complex quartic extension of $M$.
Denote by 
\[
f(x)=x^4+a_3x^3+a_2x^2+a_1x+a_0\in \Z_M[x]
\]
the relativ
defining polynomial of $\xi$ over $M$.
Assume that $K$ has a relative integral basis over $M$.
Our purpose is to determine all generators $\alpha\in\Z_K$ of relative
power integral bases of $K$ over $M$.
That is, we consider $\alpha\in \Z_K$ such that 
$(1,\alpha,\alpha^2,\alpha^3)$ is a relative integer basis of $K$
over $M$. Such an $\alpha$ has relative index
\[
I_{K/M}(\alpha)=(\Z_M[\alpha]^+:\Z_M^+)=1
\]
(see \cite{grsz}).
According to \cite{gprel4} $\alpha$ can be written in the form 
\begin{equation}
\alpha=\frac{A+X\xi+Y\xi^2+Z\xi^3}{d}
\label{alpha}
\end{equation}
where $A,X,Y,Z\in\Z_M$, and $d\in\Z$ is a non-zero common denominator.

Set
\[
F(U,V)=U^3-a_2U^2V+(a_1a_3-4a_4)UV^2+(4a_2a_4-a_3^2-a_1^2a_4)V^3,
\]
\[
Q_1(X,Y,Z)=X^2-a_1XY+a_2Y^2+(a_1^2-2a_2)XZ+(a_3-a_1a_2)YZ+(-a_1a_3+a_2^2+a_4)Z^2,
\]
\[
Q_2(X,Y,Z)=Y^2-XZ-a_1YZ+a_2Z^2.
\]
Let $i_0=I_{K/M}(\xi)$.
Recall the main result of \cite{gprel4}:

\begin{lemma}
\label{lemma1}
$\alpha\in\Z_K$ generates a relative power integral basis of $K$ over $M$
if and only if there exist $U,V\in\Z_M$ such that
\begin{equation}
N_{M/Q}(F(U,V))=\frac{d^{6m}}{i_0},
\label{F}
\end{equation}
with
\begin{equation}
Q_1(X,Y,Z)=U,\;\;\; Q_2(X,Y,Z)=V.
\label{Q12}
\end{equation}
\end{lemma}

Note that Lemma \ref{lemma1} enables us to determine finitely many
$(X,Y,Z)\in\Z_M^3$, such that
all possible generators of relative power integral bases
of $K$ over $M$ are of the form 
\begin{equation}
\alpha=\frac{A+\varepsilon(X\xi+Y\xi^2+Z\xi^3)}{d},
\label{axyz}
\end{equation}
where $\varepsilon$ is a unit in $M$ and $A\in\Z_M$ is
arbitrary (such that $\alpha\in\Z_K$).

\section{Solving the cubic relative Thue equation over $M$}
\label{cubiceq}

To apply Lemma \ref{lemma1} we first have to determine the 
solutions $U,V\in\Z_M$ of equation (\ref{F}).

By (\ref{F}) we have
\begin{equation}
F(U,V)=\varepsilon\cdot \nu
\label{fuv}
\end{equation}
where $\varepsilon$ is a unit in $M$ and $\nu\in\Z_M$ of norm $d^{6m}/i_0$.
As it is known, up to associates there are only finitely many possible
values of $\nu$ that can be determined by an algebraic number theory package like
Kash \cite{kash}, Magma \cite{magma} or Pari \cite{pari}.

\cite{gp} gives an algorithm for the resolution of relative Thue equations.
However here we would like to emphasize some special features of this calculation.
There are three possible cases according to $F$:\\
A) $F$ splits into linear factors over $M$,\\
B) $F$ is irreducible over $M$,\\
C) $F$ is a product of a linear and a quadratic factor over $M$.

We make some remarks on all the three cases, but we give a complete 
description only in the most interesting case C), with details of Baker's 
method, reduction and enumeration algorithms.

\subsection{A) $F$ splits into linear factors over $M$}
If 
\begin{equation}
F(U,V)=(U-\lambda_1 V)(U-\lambda_2 V)(U-\lambda_3 V),
\label{fuvs}
\end{equation}
with 
$\lambda_1, \lambda_2, \lambda_3\in\Z_M$, then (\ref{fuv}) implies
\begin{equation}
U-\lambda_i V=\delta_i\nu_i \;\;\; (i=1,2,3),
\label{uuvv}
\end{equation}
with units $\nu_i\in M$ and with $\delta_i\in\Z_M$ such that 
$\delta_1\delta_2\delta_3=\nu$
(the norms of $\delta_i$ divide the norm of $\nu$).
Using Siegel's identity
\[
(\lambda_1-\lambda_2)(U-\lambda_3 V)
+(\lambda_2-\lambda_3)(U-\lambda_1 V)
+(\lambda_3-\lambda_1)(U-\lambda_2 V)
=0,
\]
holding for any $U,V$.
This gives rise to a unit equation
\begin{equation}
\alpha X+\beta Y=1,
\label{unit}
\end{equation}
with
\[
\alpha=\frac{(\lambda_1-\lambda_2)\delta_3}{(\lambda_3-\lambda_2)\delta_1},\;\;
\beta=\frac{(\lambda_3-\lambda_1)\delta_2}{(\lambda_3-\lambda_2)\delta_1},
\]
where
\[
X=\frac{\nu_3}{\nu_1},\;\;
Y=\frac{\nu_2}{\nu_1}
\]
are unknown units.
This is a standard unit equation over $M$ that can be solved using the standard
methods, see \cite{book}. 
We represent $X$ and $Y$ as a power product of the fundamental units of $M$
with unknown exponents. We apply Baker's method,
reduction method and enumerate the small solutions \cite{book}.
These procedures involve only the fundamental units of $M$.

\subsection{B) $F$ is irreducible over $M$}

If $F$ is irreducible over $M$ then
it splits into linear factors over a cubic extension $L$ of $M$.
That is, (\ref{fuvs}) holds, with
$\lambda_1, \lambda_2, \lambda_3$ which are relative conjugates of
$\lambda=\lambda_1\in\Z_L$ over $M$.
(\ref{uuvv}) is valid with $\delta_i\in\Z_L$ and units $\nu_i\in L$
which are relative conjugates of $\delta,\nu\in L$ over $M$, respectively.
We obtain a unit equation (\ref{unit}) which is formally the same as above.
However there is an important difference. By calculating 
${\displaystyle
X=\frac{\nu_3}{\nu_1},\;\;
Y=\frac{\nu_2}{\nu_1}
}$
only the relative units remain in the quotients.
To explain this situation denote by $\varepsilon_1,\ldots,\varepsilon_r$ 
the fundamental units in $M$.
For the simplicity of our formulas assume that a set of fundamental 
units of $L$ is obtained by extending this set
with units $\eta_1,\ldots,\eta_s$ of $L$. That is, any unit in $L$
can be written as
\[
\varepsilon=\varepsilon_1^{a_1}\cdots \varepsilon_r^{a_r} 
\eta_1^{b_1}\cdots \eta_s^{b_s},
\]
with exponents $a_1,\ldots,a_r,b_1,\ldots,b_s$.
(Note that also in the general case any unit can be written is a similar form,
but maybe with some common denominators in the exponents, 
but from our point of view this can be dealt with analogously.)

Let $\vartheta$ be a generating element of $L$ over $M$.
Denote by $f\in\Z_M[x]$ the relative defining polynomial of $\vartheta$
over $M$. Denote by $\vartheta^{(ij)}$ ($j=1,2,3$) the roots of the $i$-th conjugate of $f$ over $M$ ($i=1,2,3$).
Denote by $\zeta^{(ij)}$ the conjugates of any $\zeta\in L$ 
corresponding to $\vartheta^{(ij)}$.
For $\delta\in M$ we have $\delta^{(ij)}=\delta^{(i)},\; j=1,2,3$.

As we have seen, in (\ref{unit}), $X$ and $Y$ contains
quotients of relative conjugates of $\nu$ over $M$, where
\[
\nu^{(ij)}=\left(\varepsilon_1^{(i)}\right)^{a_1}
\cdots
\left(\varepsilon_r^{(i)}\right)^{a_r}
\left(\eta_1^{(ij)}\right)^{b_1}\cdots\left(\eta_r^{(ij)}\right)^{b_s}.
\]
This yields that 
\[
Y=
\frac{\nu^{(ij_2)}}{\nu^{(ij_1)}}=
\left(\frac{\eta_1^{(ij_2)}}{\eta_1^{(ij_1)}}\right)^{b_1}\cdots
\left(\frac{\eta_s^{(ij_2)}}{\eta_s^{(ij_1)}}\right)^{b_s},
\]
and similarly for $X$.
Therefore we obtain a unit equation in $X$ and $Y$,
both terms with $s$ factors and the same
exponents. The standard arguments (Baker's method, reduction, enumeration) can be
used to calculate $b_1,\ldots,b_s$. This can be used to determine $U,V$ up to a unit factor in $M$.
One can apply Baker's method, reduction and enumeration the same standard way as
in A) but with the above $s$ factors corresponding to the relative units.

\subsection{C) $F$ is a product of a linear and a quadratic factor over $M$}

The most interesting case is, when
$F$ is a product of a linear and a quadratic factor over $M$.
Then we have
\[
F(U,V)=(U-\lambda_1 V)(U^2+\lambda_2 UV+\lambda_3 V^2),
\]
with $\lambda_i\in \Z_M\; (i=1,2,3)$ where the second degree factor
is irreducible over $M$. Denote by $G=M(\gamma)$ ($\gamma\in\Z_G$)
a quadratic extension of $M$ such 
that the quadratic factor of $F$ splits into linear factors over $G$.
Denote by $\gamma^{(ij)}$ ($j=1,2$) the roots of the $i$-th conjugate
of the relative defining polynomial of $\gamma$ over $M$ ($i=1,\ldots,m$).
Denote by $\delta^{(ij)}$ the conjugates of any $\delta\in G$ corresponding
to $\gamma^{(ij)}$ ($i=1,\ldots,m,j=1,2$).
For $\zeta\in M$ we have $\zeta^{(ij)}=\zeta^{(i)},j=1,2$.
Then we have 
\begin{equation}
F^{(ij)}(U,V)=(U-\lambda_1^{(i)} V)(U-\gamma^{(i1)}V)(U-\gamma^{(i2)}V).
\label{Fi}
\end{equation}
By (\ref{F}) and (\ref{Fi}) we have
\begin{eqnarray}
U-\lambda_1^{(i)} V&=& \delta_M^{(i)}\nu_M^{(i)}  \nonumber\\
U-\gamma^{(i1)}V&=&\delta_G^{(i1)}\nu_G^{(i1)}\label{xy}\\
U-\gamma^{(i2)}V&=&\delta_G^{(i2)}\nu_G^{(i2)}\nonumber
\end{eqnarray}
where $\delta_M\in \Z_M$, the norm of which devides $d^{6m}/i_0$,
$\nu_M$ is a unit in $M$, 
$\delta_G\in \Z_G$, the norm of which devides $d^{6m}/i_0$
and $\nu_G$ is a unit in $G$. 
(Up to associates there are only a few possible values of $\delta_M,\delta_G$.)
Siegel's identity gives
\[
(\lambda_1^{(i)}-\gamma^{(i1)})(U-\gamma^{(i2)}V)
+(\gamma^{(i1)}-\gamma^{(i2)})(U-\lambda_1^{(i)}V)
+(\gamma^{(i2)}-\lambda_1^{(i)})(U-\gamma^{(i1)}V)
=0,
\]
whence
\begin{equation}
\alpha X+\beta Y=1,
\label{abxy}
\end{equation}
with
\[
\alpha=\frac{(\lambda_1^{(i)}-\gamma^{(i1)})\delta_G^{(i2)}}
{(\gamma^{(i2)}-\gamma^{(i1)})\delta_M^{(i)}},
\;\;
\beta=\frac{(\gamma^{(i2)}-\lambda_1^{(i)})\delta_G^{(i1)}}
{(\gamma^{(i2)}-\gamma^{(i1)})\delta_M^{(i)}},
\]
and the unknown units are
\[
X=\frac{\nu_G^{(i2)}}{\nu_M^{(i)}},
\;\;
Y=\frac{\nu_G^{(i1)}}{\nu_M^{(i)}}.
\]
Observe that $X$ and $Y$ are conjugated over $M$, as well as
$\alpha$ and $\beta$.
Let again $\varepsilon_1,\ldots\varepsilon_r$ denote the fundamental units of 
$M$. For simplicity's sake 
assume that a system of fundamental units of $G$
is obtained by extending this system by some relative units
$\eta_1,\ldots,\eta_s$.
Set
\[
X^{(ij)}=(\varepsilon_1^{(i)})^{a_1}\cdots (\varepsilon_r^{(i)})^{a_r}\cdot
(\eta_1^{(ij)})^{b_1}\cdots (\eta_s^{(ij)})^{b_s}.
\]
Let $A=\max |a_i|, B=\max |b_i|, E=\max(A,B)$.

We detail the application of Baker's method, reduction and enumeration
algorithms in this case C).

\subsubsection{Baker's method}

We apply standard arguments.
Since $|N_{G/Q}(X)|=1$, hence $\sum_{i=1}^m\sum_{j=1}^2 |\log |X^{ij}||=0$.
There exists a conjugate $i_0,j_0$ with
\[
\log |X^{(i_0,j_0)}|<- c_1 E,
\]
where $c_1$ is a positive constant that can be easily calculated.
According to standard arguments (cf. \cite{book}) we have
\[
\exp(-c_1 E)> |X^{(i_0,j_0)}|=
\frac{1}{|\alpha^{(i_0,j_0)}|}|1-\beta^{(i_0,j_0)}Y^{(i_0,j_0)}|
\geq 
\frac{1}{2|\alpha^{(i_0,j_0)}|} |\log|\beta^{(i_0,j_0)}Y^{(i_0,j_0)}||
\]
\begin{equation}
\geq
\frac{1}{2|\alpha^{(i_0,j_0)}|}
|\log |\beta^{(i_0,j_0)}|
+a_1\log |\varepsilon_1^{(i_0)}|+\ldots +a_r\log |\varepsilon_r^{(i_0)}|
+b_1\log |\eta_1^{(i_0,j_0)}|+\ldots +b_s\log |\eta_s^{(i_0,j_0)}||
\label{i2}
\end{equation}
\[
\geq
\exp (-C \log E),
\]
where $C$ is a huge positive constant.
In the last step we applied Baker's method, that is e.g. the 
estimates of A. Baker and G. W\"ustholz \cite{bawu}.
Comparing the beginning and the end of this series of
inequalities we obtain an upper bound for $E$.
$E_B$ is going to be the maximum of these upper bounds 
(the maximum taken for all possible pairs ($i_0,j_0$)).
This usually of magnitude $10^{30}-10^{100}$.

\subsubsection{Reduction}

The next step is to reduce the bound $E_B$. For this purpose we apply
Lemma 2.2.2 of \cite{book}. We recall here the cruxial statement.

Let $\zeta_1,\ldots, \zeta_n$ be multiplicative independent
algebraic numbers, $d_1,\ldots,d_n$ integers. Set $D=\max |d_i|$.
Assume that 
\begin{equation}
|d_1\zeta_1+\ldots+d_n\zeta_n|<c_1\exp(-c_2 D-c_3),
\label{redineq}
\end{equation}
for any $d_1,\ldots,d_n$ where $c_1,c_2,c_3$ 
are given positive constants (of moderate size).
Our purpose is to reduce the bound $D_0$ obtained
previously for $D$ by using Baker's method.

Let $H$ be a large constant (an appropriate value is about $D_0^n$) and consider the
lattice $\cal L$ spanned by the columns of the $n+2$ by $n$ matrix
\[
\left(
\begin{array}{cccc}
1&0&\ldots &0\\
0&1&\ldots &0\\
\vdots &\vdots &\vdots &\vdots \\
0&0&\ldots&1\\
H\cdot{\rm Re}(\zeta_1)&H\cdot{\rm Re}(\zeta_2)&\ldots &H\cdot{\rm Re}(\zeta_n)\\
H\cdot{\rm Im}(\zeta_1)&H\cdot{\rm Im}(\zeta_2)&\ldots &H\cdot{\rm Im}(\zeta_n)\\
\end{array}
\right).
\]
Assume that the columns in the above matrix are linearly independent.
Denote by $b_1$ the first vector of an LLL-reduced basis of this lattice
(cf. A. K. Lenstra, H. W. Lenstra Jr. and L. Lov\'asz \cite{lll},
M. Pohst \cite{dmv}).

\begin{lemma}
\label{redlemma}
If $D\leq D_0$ and $|b_1|\geq \sqrt{(n+1)2^{n-1}} \cdot D_0$, then
\[
D\leq \frac{\log H+\log c_1-c_3-\log D_0}{c_2}.
\]
\end{lemma}

\vspace{0.5cm}

We apply this lemma based on the inequality
\[
|\log |\beta^{(i_0,j_0)}|
+a_1\log |\varepsilon_1^{(i_0)}|+\ldots +a_r\log |\varepsilon_r^{(i_0)}|
+b_1\log |\eta_1^{(i_0,j_0)}|+\ldots +b_s\log |\eta_s^{(i_0,j_0)}||
\]
\[
<2|\alpha^{(i_0,j_0)}|\cdot \exp(-c_1 E),
\]
see (\ref{i2}), letting
$n=r+s+1, \zeta_1=\log|\beta^{(i_0,j_0)}|,
\zeta_2=\log |\varepsilon_1^{(i_0)}|,\ldots,
\zeta_{r+1}=\log |\varepsilon_r^{(i_0)}|$,
$\zeta_{r+2}=\log |\eta_1^{(i_0,j_0)}|,\ldots,
\zeta_{r+s+1}=\log |\eta_r^{(i_0,j_0)}|$.
In the totally real case we may omit the last row of the matrix.

The repeated application of the reduction method brings down
the bound obtained by Baker's method to a so called reduced bound.
We perform this calculation for all possible pairs ($i_0,j_0$)
and $E_R$ is the maximum of the reduced bounds. This is usually
of magnitude 100-1000 depending on the size of the example.

\subsubsection{Enumeration}

Observe that although the reduced bound $E_R$ is rather small,
the number of possible values $-E_R\leq a_1,\ldots,a_r,b_1,\ldots b_r
\leq E_R$ is
a huge number: $(2E_R+1)^{r+s}$.

Hence we must apply the enumeration methods of \cite{book}.

As we shall see, in the case C) that we consder,
the general method of \cite{book} can only be applied with 
non-trivial modifications.

Recall that in our equation (\ref{abxy})
$X$ and $Y$ are conjugated over $M$, as well as
$\alpha$ and $\beta$.
We assumed that 
\[
X^{(ij)}=(\varepsilon_1^{(i)})^{a_1}\cdots (\varepsilon_r^{(i)})^{a_r}\cdot
(\eta_1^{(ij)})^{b_1}\cdots (\eta_s^{(ij)})^{b_s},
\]
where $\varepsilon_1,\ldots,\varepsilon_r$ is a set of fundamental units of
$M$ and 
$\varepsilon_1,\ldots,\varepsilon_r,\eta_1\ldots,\eta_s$ is a set of
fundamental units in $G$.
Further $A=\max |a_i|, B=\max |b_i|, E=\max(A,B)$.

We write equation (\ref{abxy}) in the form
\begin{equation}
\alpha^{(i1)}X^{(i1)}+\alpha^{(i2)}X^{(i2)}=1.
\label{axax}
\end{equation}
The essence of the enumeration algorithm of \cite{book} is that we calculate
an $S$ with
\begin{equation}
\frac{1}{S}\leq |\alpha^{(ij)} X^{(ij)}|\leq S,
\label{SS}
\end{equation}
for all conjugates. 
As an initial value we can take
\[
\log S=\max_{i,j}|\alpha^{(ij)}|+
E_R\cdot \max_{i,j}
(\log|\varepsilon_1^{(i)}|+\cdots +\log|\varepsilon_r^{(i)}|
+\log|\eta_1^{(ij)}|+\cdots +\log|\eta_s^{(ij)}|).
\]
Then we show that $S$ can be replaced by a smaller
constant $s$ on the prize that we enumerate and test the exponent vectors
$(a_1,\ldots,a_r,b_1,\ldots,b_s)$ 
in some exceptional sets. These sets are ellipsoids
containing relatively few possible vectors. We repeat diminishing $S$ to $s$
until we reach a relatively small $s$ value (of magnitude 10).
Finally we enumerate an ellipsoid of type (\ref{SS}) with the final $S=s$.

We describe now how this procedure can be adopted to our case. 
Consider the following Lemma which is a modified version of Lemma 2.3.1
of \cite{book}.

\begin{lemma}
Let $s<S$ and assume that (\ref{SS}) holds. If there exists a $j$
such that 
\[
\frac{1}{s}\leq |\alpha^{(ij)}X^{(ij)}|\leq s
\]
is violated, then \\
I. either there exists a $j_0$ with
\begin{equation}
|\log|\alpha^{(ij_0)}X^{(ij_0)}||\leq \frac{2}{s},
\label{s1}
\end{equation}
II. or there exists a $j_0$ such that 
\begin{equation}
\left|\log\left|\frac{\alpha^{(ij_0)}X^{(ij_0)}}{\alpha^{(ij)}X^{(ij)}}\right|\right|
\leq \frac{2}{s}.
\label{s2}
\end{equation}
\label{enumlemma}
\end{lemma}

\noindent
{\bf Proof}\\
Let $j_0=\{1,2\}\setminus \{j\}$.\\
If 
\[
\frac{1}{S}\leq |\alpha^{(ij)}X^{(ij)}|\leq \frac{1}{s},
\]
then we have
\[
|\log|\alpha^{(ij_0)}X^{(ij_0)}||\leq 2||\alpha^{(ij_0)}X^{(ij_0)}|-1|
\leq 2|\alpha^{(ij_0)}X^{(ij_0)}-1|=2|\alpha^{(ij)}X^{(ij)}|
\leq \frac{2}{s},
\]
which implies (\ref{s1}).

On the other hand, if 
\[
s\leq |\alpha^{(ij)}X^{(ij)}|\leq S,
\]
then
\[
\left|\log\left|\frac{\alpha^{(ij_0)}X^{(ij_0)}}{\alpha^{(ij)}X^{(ij)}}\right|\right|
\leq 2\left|1-\left|\frac{\alpha^{(ij_0)}X^{(ij_0)}}{\alpha^{(ij)}X^{(ij)}}\right|\right|
\leq 2\left|1+\frac{\alpha^{(ij_0)}X^{(ij_0)}}{\alpha^{(ij)}X^{(ij)}}\right|
=\frac{2}{|\alpha^{(ij)}X^{(ij)}|}\leq \frac{2}{s},
\]
which implies (\ref{s2}).

Remark that in both cases we used the inequality $|\log x|<2|x-1|$
holding for any complex number $x$ with $|x-1|<0.795$. In our calculations 
this is satisfied, since in the applications we have $2/s<0.795$.
$\Box$

Now we explain how to enumerate the possible vectors $a_1,\ldots,a_r,
b_1,\ldots,b_s$ if in addition to (\ref{SS}) (for all $i,j$) either 
(\ref{s1}) or (\ref{s2}) is satisfied (for certain $i,j_0$).

\vspace{1cm}

\noindent
{\bf Case I.}\\
We have
\[
\log|\alpha^{(ij)}X^{(ij)}|
=\log|\alpha^{(ij)}|+a_1\log |\varepsilon_1^{(i)}|
+\ldots +a_r\log |\varepsilon_r^{(i)}|
+b_1\log|\eta_1^{(ij)}|+\ldots +b_s \log|\eta_s^{(ij)}|,
\]
for all $i,j$.
Let
\[
h:=\left(
\begin{array}{c}
\log|\alpha^{(11)}X^{(11)}|\\
\log|\alpha^{(12)}X^{(12)}|\\
\log|\alpha^{(21)}X^{(21)}|\\
\log|\alpha^{(22)}X^{(22)}|\\
\vdots\\
\log|\alpha^{(m1)}X^{(m1)}|\\
\log|\alpha^{(m2)}X^{(m2)}|\\
\log|\alpha^{(ij_0)}X^{(ij_0)}|
\end{array}
\right)
,\;\;
g=
\left(
\begin{array}{c}
\log|\alpha^{(11)}|\\
\log|\alpha^{(12)}|\\
\log|\alpha^{(21)}|\\
\log|\alpha^{(22)}|\\
\vdots\\
\log|\alpha^{(m1)}|\\
\log|\alpha^{(m2)}|\\
\log|\alpha^{(ij_0)}|
\end{array}
\right)
\]
\[
e_k=
\left(
\begin{array}{c}
\log|\varepsilon^{(11)}|\\
\log|\varepsilon^{(12)}|\\
\log|\varepsilon^{(21)}|\\
\log|\varepsilon^{(22)}|\\
\vdots\\
\log|\varepsilon^{(m1)}|\\
\log|\varepsilon^{(m2)}|\\
\log|\varepsilon^{(ij_0)}|
\end{array}
\right)
,\; (k=1,\ldots,r),
\;\;\;
f_l=
\left(
\begin{array}{c}
\log|\eta^{(11)}|\\
\log|\eta^{(12)}|\\
\log|\eta^{(21)}|\\
\log|\eta^{(22)}|\\
\vdots\\
\log|\eta^{(m1)}|\\
\log|\eta^{(m2)}|\\
\log|\eta^{(ij_0)}|
\end{array}
\right)
,\; (l=1,\ldots,s)
\]
Then
\[
h=g+a_1e_1+\ldots a_re_r+b_1f_1+\ldots b_sf_s
\]
Let
\[
\lambda_{ij}=\frac{1}{\log S}, (1\leq i\leq m,j=1,2),
\]
and for the last coordinate let
\[
\lambda=\frac{s}{2}.
\]
For any vector
\[
v=\left(
\begin{array}{c}
x_{11}\\
x_{12}\\
x_{21}\\
x_{22}\\
\vdots\\
x_{m1}\\
x_{m2}\\
x_{ij_0}
\end{array}
\right)
\;\;\; {\rm set} \;\;\; 
\varphi (v)=\left(
\begin{array}{c}
\lambda_{11}\cdot x_{11}\\
\lambda_{12}\cdot x_{12}\\
\lambda_{21}\cdot x_{21}\\
\lambda_{22}\cdot x_{22}\\
\vdots\\
\lambda_{m1}\cdot x_{m1}\\
\lambda_{m2}\cdot x_{m2}\\
\lambda \cdot x_{ij_0}
\end{array}
\right).
\]
Then
\[
\varphi(h)=\varphi(g)+a_1\varphi(e_1)+\ldots a_r\varphi(e_r)
+b_1\varphi(f_1)+\ldots b_s\varphi(f_s).
\]
Further, (\ref{SS}) and (\ref{s1}) imply
\[
||\varphi(h)||^2=
||\varphi(g)+a_1\varphi(e_1)+\ldots a_r\varphi(e_r)
+b_1\varphi(f_1)+\ldots b_s\varphi(f_s)
||^2
\]
\[
\leq \sum_{i=1}^m\sum_{j=1}^2\frac{1}{\log S}|\log|\alpha^{(ij)}X^{(ij)}||
+\frac{s}{2}|\log|\alpha^{(ij_0)}X^{(ij_0)}||
\leq 2m+1.
\]
The above $L^2$ norm defines an ellipsoid.

\vspace{1cm}

\noindent
{\bf Case II.}\\
We have
\[
\log\left|\frac{\alpha^{(ij_0)}X^{(ij_0)}}{\alpha^{(ij)}X^{(ij)}}\right|
=
\log\left|\frac{\alpha^{(ij_0)}}{\alpha^{(ij)}}\right|
+
b_1\log\left|\frac{\eta_1^{(ij_0)}}{\eta_1^{(ij)}}\right|
+\ldots +
b_s\log\left|\frac{\eta_s^{(ij_0)}}{\eta_s^{(ij)}}\right|.
\]
Observe that here we only have quotients of conjugates
of the relative units $\eta_1,\ldots,\eta_s$.
By (\ref{SS}) we can derive 
\begin{equation}
\frac{1}{S^2}<
\left|\frac{\alpha^{(ij_0)}X^{(ij_0)}}{\alpha^{(ij)}X^{(ij)}}\right|
<S^2
\label{S2}
\end{equation}
for any $i,j$. Further, (\ref{s2}) holds for $i,j_0$.
Let
\[
h:=\left(
\begin{array}{c}
\log\left|\frac{\alpha^{(11)}X^{(11)}}{\alpha^{(12)}X^{(12)}}\right|\\ \\
\log\left|\frac{\alpha^{(21)}X^{(21)}}{\alpha^{(22)}X^{(22)}}\right|\\
\vdots\\
\log\left|\frac{\alpha^{(m1)}X^{(m1)}}{\alpha^{(m2)}X^{(m2)}}\right|\\ \\
\log\left|\frac{\alpha^{(ij_0)}X^{(ij_0)}}{\alpha^{(ij)}X^{(ij)}}\right|
\end{array}
\right)
,\;\;
g=
\left(
\begin{array}{c}
\log\left|\frac{\alpha^{(11)}}{\alpha^{(12)}}\right|\\ \\
\log\left|\frac{\alpha^{(21)}}{\alpha^{(22)}}\right|\\
\vdots\\
\log\left|\frac{\alpha^{(m1)}}{\alpha^{(m2)}}\right|\\ \\
\log\left|\frac{\alpha^{(ij_0)}}{\alpha^{(ij)}}\right|
\end{array}
\right)
,\;\;
f_l=
\left(
\begin{array}{c}
\log\left|\frac{\varepsilon_l^{(11)}}{\varepsilon_l^{(12)}}\right|\\ \\
\log\left|\frac{\varepsilon_l^{(21)}}{\varepsilon_l^{(22)}}\right|\\
\vdots\\
\log\left|\frac{\varepsilon_l^{(m1)}}{\varepsilon_l^{(m2)}}\right|\\ \\
\log\left|\frac{\varepsilon_l^{(ij_0)}}{\varepsilon_l^{(ij)}}\right|
\end{array}
\right)
,\; (l=1,\ldots,s).
\]
Then
\[
h=g+b_1f_1+\ldots +b_sf_s
\]
Let
\[
\lambda_{i}=\frac{1}{2\log S}, (1\leq i\leq m),
\]
and for the last coordinate let
\[
\lambda_{m+1}=\frac{s}{2}.
\]
For any vector
\[
v=\left(
\begin{array}{c}
x_{1}\\
x_{2}\\
\vdots\\
x_{m}\\
x_{m+1}
\end{array}
\right)
\;\;\; {\rm set} \;\;\; 
\varphi (v)=\left(
\begin{array}{c}
\lambda_1\cdot x_{1}\\
\lambda_2\cdot x_{2}\\
\vdots\\
\lambda_m\cdot x_{m}\\
\lambda_{m+1}\cdot x_{m+1}
\end{array}
\right).
\]
Then
\[
\varphi(h)=\varphi(g)+b_1\varphi(f_1)+\ldots b_s\varphi(f_s).
\]
Further, (\ref{S2}) and (\ref{s2}) imply
\[
||\varphi(h)||^2=
||\varphi(g)+b_1\varphi(f_1)+\ldots b_s\varphi(f_s)||^2
\]
\[
\leq \sum_{i=1}^m \frac{1}{2\log S}
\left|\log\left|\frac{\alpha^{(i1)}X^{(i1)}}{\alpha^{(i2)}X^{(i2)}}\right|\right|
+\frac{s}{2}
\left|\log\left|\frac{\alpha^{(ij_0)}X^{(ij_0)}}{\alpha^{(ij)}X^{(ij)}}\right|\right|
\leq m+1.
\]
The above $L^2$ norm defines an ellipsoid.

\vspace{1cm}

{\bf Remarks}

\begin{enumerate}
\item
The procedure can be continued both in Case I and Case II by taking 
the previous $s$ in the role of $S$ and choosing a smaller $s$.
For an appropriate choice of $s$ see \cite{book}. Usually
we take $s=\sqrt{S}$.

\item
Proceeding until a relatively small value of $S$ (of magnitude 10),
we enumerate in both cases an ellipsoid taking all weights $1/\log S$ in Case I
and $1/(2\log S)$ in Case II to finish the procedure (cf. \cite{book},
see the last ellipsoids in our example).

\item
Enumerating the ellipsoids, in Case I we obtain all possible 
exponent vectors $(a_1,\ldots,a_r, b_1,\ldots,b_s)$.
Observe that in Case II we can only enumerate the possible 
values of $(b_1,\ldots,b_s)$.
For all possible $(b_1,\ldots,b_s)$ we run $a_1,\ldots,a_r$
run between $-E_R$ and $E_R$ and test if the unit equation (\ref{axax}) holds.
Since the ground field $M$ is usually of small degree, 
this can be done relatively fast.

\item
We emphasize that for the enumeration of
the ellipsoids we use the improved method involving LLL reduction
(see \cite{dmv}).

\item
To speed up the test of possible exponent vectors we use sieves
(see \cite{book}). This enables us to check mod $p$ congruences instead of
equations in high precision real numbers. Also, calculations with integers modulo $p$
is much faster than real arithmetic.

\item
The enumeration of exponent vectors in the exceptional ellipsoids must be
performed for all possible pairs ($i,j_0$) and we have to check all possible
exponent vectors.

\end{enumerate}

\section{The quartic relative Thue equation}
\label{qqeq}

Having $U, V$ (determined up to a unit factor in $M$)
we follow the methods of \cite{gprel4} 
(see also \cite{book}) to determine $X,Y,Z$.
Set 
\[
Q_0(X,Y,Z)=U\cdot Q_1(X,Y,Z)-V\cdot Q_2(X,Y,Z).
\]
Let $X_0,Y_0,Z_0\in\Z_M$ be a nontrivial solution of 
\begin{equation}
Q_0(X,Y,Z)=0,
\label{q00}
\end{equation} 
with, say $Z_0\ne 0$.
We represent $X,Y,Z$ with parameters $P,Q,R\in M$
in the form
\begin{eqnarray}
X&=&R\cdot X_0+P,\nonumber\\
Y&=&R\cdot Y_0+Q,\label{rpq}\\
Z&=&R\cdot Z_0.\nonumber
\end{eqnarray}
Substituting these representations into (\ref{q00})
we obtain an equation of the form
\[
R(C_1P+C_2Q)=C_3P^2+C_4PQ+C_5Q^2,
\]
with $C_1,\ldots,C_5\in\Z_M$. We multiply the equations (\ref{rpq}) 
by $C_1P+C_2Q$ and use the above equation to eliminate $R$
on the right hand sides. We obtain 
\begin{eqnarray}
\kappa\cdot X &=& f_X(P,Q),\nonumber \\
\kappa\cdot Y &=& f_Y(P,Q),\label{xyzpq}\\
\kappa\cdot Z &=& f_Z(P,Q),\nonumber
\end{eqnarray}
with quadratic form $f_X,f_Y,f_Z\in\Z_M[P,Q]$.
As stated in \cite{gprel4} we can replace $P,Q,\kappa$ by integer parameters
(by multiplying the equations by the square of a common denominator 
of $\kappa,P,Q$) 
and $\kappa$ may attain only finitely many non-associated values.
Substituting these representations into (\ref{Q12}) we obtain
quartic equations over $M$:
\begin{eqnarray}
F_1(P,Q)=Q_1(f_X(P,Q),f_Y(P,Q),f_Z(P,Q))&=&\kappa^2\cdot U,\nonumber\\
F_2(P,Q)=Q_2(f_X(P,Q),f_Y(P,Q),f_Z(P,Q))&=&\kappa^2\cdot V\label{relthue}.
\end{eqnarray}
According to \cite{gprel4} 
at least one of these is a quartic relative Thue equation
over $M$, having a root in $K$.

\section{Relative Thue equations in totally complex extensions
of totally real fields}

We show that to solve relative Thue equations of type (\ref{relthue})
is a trivial matter. We formulate our assertion in a general form. 

Let $M$ be a totally real number field and let $K=M(\xi)$ be a 
totally complex extension of degree $k$ of $M$. More exactly, if $f(x)\in\Z_M[x]$
is the relative defining polynomial of $\xi$ over $M$, then
all roots $\xi^{(i1)},\ldots, \xi^{(ik)}$ of a conjugate $f^{(i)}(x)$
of $f(x)$ are complex. Set 
\[
c_0=\frac{1}{\min_{i,j}|{\rm Im}(\xi^{(ij)})|}.
\]
Let $0\ne \nu\in\Z_M$, set
\[
F(X,Y)=N_{K/M}(X-\xi Y),
\]
and consider the relative Thue equation
\begin{equation}
F(X,Y)=\nu, \;\; {\rm in}\;\; X,Y\in\Z_M.
\label{ff}
\end{equation}
We denote by $|\overline{Z}|$ the size of $Z\in M$, that is the
maximum absolute value of its conjugates.

\begin{theorem}
All solutions $X,Y\in\Z_M$ of equation (\ref{ff}) satisfy
\begin{equation}
\max (|\overline{X}|,|\overline{Y}|)
\leq |\overline{\nu}|^{1/k}(1+c_0 |\overline{\xi}|).
\label{xysize}
\end{equation}
\label{thth}
\end{theorem}

\noindent
{\bf Proof.}\\
Let $X,Y\in\Z_M$ be an arbitrary but fixed solution of equation (\ref{ff}).
Denote by $\gamma^{(i)}$ the conjugates of $\gamma\in M$ corresponding to 
$\xi^{(i1)},\ldots, \xi^{(ik)}$.
For any $1\leq i\leq m,\; 1\leq j\leq k$ 
set
\[
\beta^{(ij)}=X^{(i)}-\xi^{(ij)}Y^{(i)}.
\]
For any $1\leq i\leq m$ we have
\[
\prod_{j=1}^k \beta^{(ij)}=\nu^{(i)}.
\]
Let $j_0$ be the index with
\[
|\beta^{(ij_0)}|=\min_{1\leq j\leq k}|\beta^{(ij)}|.
\]
Then 
\[
|{\rm Im}(\xi^{(ij_0)})|\cdot |Y^{(i)}|
={\rm Im}|X^{(i)}-\xi^{(ij_0)}Y^{(i)}|
={\rm Im}|\beta^{(ij_0)}|
\leq |\nu^{(i)}|^{1/k},
\]
whence
\[
|Y^{(i)}|\leq \frac{|\nu^{(i)}|^{1/k}}{|{\rm Im}(\xi^{(ij_0)})|},\;\;
|X^{(i)}|\leq |\nu^{(i)}|^{1/k}+|\xi^{(ij_0)}|\cdot |Y^{(i)}|,
\]
which implies our assertion.\\
$\Box$

\vspace{1cm}

Let us return now to the equations (\ref{relthue}).
We denoted the fundamental units of $M$ by $\varepsilon_1,\ldots,\varepsilon_r$.
$\kappa$ can be written in the form
\[
\kappa=\kappa_0\cdot \varepsilon_1^{k_1}\ldots \varepsilon_r^{k_r}
\] 
where $\kappa_0$ can only take finitely many values, $k_1,\ldots,k_r\in\Z$.
Set $k_i=4k_i'+\ell_i$ with $-1\leq \ell_i\leq 2,\; (i=1,\ldots,r)$ and
\[
P'=P\cdot \varepsilon_1^{-k_1'}\cdots \varepsilon_r^{-k_r'},
\;\;
Q'=Q\cdot \varepsilon_1^{-k_1'}\cdots \varepsilon_r^{-k_r'}.
\]
Then we have 
\begin{eqnarray*}
F_1(P',Q')&=&\kappa_0\cdot U\cdot \varepsilon_1^{\ell_1}\cdots \varepsilon_r^{\ell_r},
\\
F_2(P',Q')&=&\kappa_0\cdot V\cdot \varepsilon_1^{\ell_1}\cdots \varepsilon_r^{\ell_r}.
\end{eqnarray*}
One of these equations is a quartic relative Thue equation (see \cite{book})
that can be solved easily by above Theorem for all possible values of $\kappa_0$
and for all possible $\ell_1,\ldots,\ell_r$. This gives $P$ and $Q$
up to a unit factor of $M$, whence we obtain $X,Y,Z$ by (\ref{xyzpq})
up to a unit factor in $M$.
The generators of relative power integral bases of $K$ over $M$
are obtained by (\ref{axyz}).
The possible values of $\alpha$ must be checked.

\section{An example}

As an example consider a root $\tau$ of the polynomial 
\[
f_3(x)=x^3-2x^2-5x-1.
\]
This is a totally real cubic field (one of the "simplest cubic fields" of
D. Shanks \cite{shanks}). The conjugates of $\mu=\tau+2$ are all positive,
$\mu$ having defining polynomial $g_3(x)=x^3-8x^2+15x-7$.
Therefore $\xi=\sqrt[4]{-\mu}$ is a totally complex algebraic integer
of degree 4 over $M$.
Note that $\xi$ has defining polynomial 
\[
f_{12}(x)=x^{12}+8x^8+15x^4+7,
\]
and the field $K=\Q(\xi)$ has an integral basis 
$(1,\xi,\ldots,\xi^{11})$, a power integral basis.
Hence $K$ is monogenic.
Our purpose is to determine all (non-equivalent) generators of
relative power integral bases of $K$ over $M$.

The integral basis $(1,\xi,\ldots,\xi^{11})$ of $K$
implies that the discriminant of $K$ is 
\[
D_K=D(f_{12})=2^{24}\cdot 7^3\cdot 19^8.
\]
It is easy to check, that
\[
(1,\mu,\mu^2,\xi,\xi\mu,\xi\mu^2,
  \xi^2,\xi^2\mu,\xi^2\mu^2, 
\xi^3,\xi^3\mu,\xi^3\mu^2)
\]
is also an integral basis of $K$, therefore any $\alpha\in\Z_K$
can be written in the form
\begin{equation}
\alpha=A+X\xi+Y\xi^2+Z\xi^3,
\label{AAlpha}
\end{equation}
with $A,X,Y,Z\in\Z_M$ (that is the constant $d$ in Lemma \ref{lemma1} is 1).
Moreover, $m=3$ and  $i_0=I_{K/M}(\xi)=1$. 
The relative defining polynomial 
of $\xi$ over $M$ is $x^4+\mu$, therefore equation (\ref{F}) is of the form
\[
N_{M/\Q}(U(U^2-4\mu V^2))=\pm 1,
\]
that is
\begin{equation}
U(U^2-4\mu V^2)=\nu,
\label{Fex}
\end{equation}
where $\nu$ is a unit in $M$.

In the absolute case a similar equation is trivial to solve, that is not
the case here. Actually the above equation is of Case C) of 
the three main types of equations detailed in Section \ref{cubiceq}.

Let $\gamma=\sqrt{\mu}$, $G=M(\gamma)$. This is a totally real sextic number field, and our
$F(U,V)$ factorizes over $\Z_G$:
\[
U(U+2\gamma V)(U-2\gamma V)=\nu.
\]
This implies, that all three factors are units, the second and third factors
are units in $G$, conjugated over $M$. 
For $i=1,2,3$ we obtain
\begin{eqnarray}
U^{(i)}&=&\nu_M^{(i)},\nonumber\\
U+2\gamma^{(i1)}V&=&\nu_G^{(i1)},\label{ccc}\\
U-2\gamma^{(i1)}V&=&U+2\gamma^{(i2)}V=\nu_G^{(i2)}.\nonumber
\end{eqnarray}
where
$\nu_M$ is a unit in $M$, $\nu_G$ is a unit in $G$.
By the above equations we obtain
\[
2\nu_M^{(i)}=\nu_G^{(i1)}+\nu_G^{(i2)}
\]
that is 
\begin{equation}
\frac{1}{2} X^{(i1)}+\frac{1}{2} X^{(i2)}=1,
\label{uuu}
\end{equation}
$X=\nu_G/\nu_M$ is a unit in $G$.

Using Kash \cite{kash} we calculated a system of fundamental units in $G$.
$\gamma=\sqrt{\mu}$ has defining polynomial $f_6(x)=x^6-8x^4+15x^2-7$.
The elements $(1,\gamma,\gamma^2,\gamma^3,\gamma^4,\gamma^5)$
form an integral basis in $G$. The coefficients of the fundamental units
in this integral basis are
\begin{eqnarray*}
\varepsilon_1:&& \;\; (1,1,0,0,0,0),\\
\varepsilon_2:&& \;\; (1,-1,0,0,0,0),\\
\varepsilon_3:&& \;\; (2,0,-1,0,0,0),\\
\varepsilon_4:&& \;\; (3,1,-1,0,0,0),\\
\varepsilon_5:&& \;\; (3,-7,0,-7,0,-1).
\end{eqnarray*}
The relative conjugates (over $M$) of these units are the following:
\begin{eqnarray*}
\varepsilon_1^{(i2)}&=&\varepsilon_2^{(i1)},\\
\varepsilon_2^{(i2)}&=&\varepsilon_1^{(i1)},\\
\varepsilon_3^{(i2)}&=&\varepsilon_3^{(i1)},\\
\varepsilon_4^{(i2)}&=&
\frac{\varepsilon_3^{(i1)}}
{\varepsilon_1^{(i1)} \varepsilon_2^{(i1)} \varepsilon_4^{(i1)}},\\
\varepsilon_5^{(i2)}&=&
\frac{\varepsilon_3^{(i1)}}
{\varepsilon_1^{(i1)} \varepsilon_2^{(i1)} \varepsilon_5^{(i1)}}.
\end{eqnarray*}
Equation (\ref{uuu}) can be written as
\begin{eqnarray}
&&
\pm
\frac{1}{2}\cdot
(\varepsilon_1^{(i1)})^{a_1}\cdot
(\varepsilon_2^{(i1)})^{a_2}\cdot
(\varepsilon_3^{(i1)})^{a_3}\cdot
(\varepsilon_4^{(i1)})^{a_4}\cdot
(\varepsilon_5^{(i1)})^{a_5}
 \nonumber \\&&
\pm
\frac{1}{2}\cdot
(\varepsilon_1^{(i2)})^{a_1}\cdot
(\varepsilon_2^{(i2)})^{a_2}\cdot
(\varepsilon_3^{(i2)})^{a_3}\cdot
(\varepsilon_4^{(i2)})^{a_4}\cdot
(\varepsilon_5^{(i2)})^{a_5}
=1
\label{u2}
\end{eqnarray}

In our example we have $c_1=0.18298,\;\; C=2.85992\cdot 10^{28}$
(cf. (\ref{i2})).
Comparing the upper and lower bounds of the above series of inequalities
(\ref{i2}) we obtain $A_B<10^{32}$.

The following table summarizes the reduction procedure:
\[
\begin{array}{|c||c|c|c|c|}\hline
                 & A<      &  H=      &  {\rm precision}  & {\rm new \; bound \; for}\; A \\ \hline
{\rm Step \;\;I} & 10^{32} & 10^{170} & 300\;{\rm digits}  &  1736 \\ \hline
{\rm Step \;\;II}& 1736    & 10^{30}  & 100\;{\rm digits}  &  336  \\ \hline
{\rm Step \;\;III}& 336    & 10^{20}  & 100\;{\rm digits}  &  219\\ \hline
\end{array}
\]
The procedure took a few minutes all together and we obtained the reduced bound
$A_R=219$.
Observe that although the reduced bound is rather small,
the number of possible values $-219\leq a_1,\ldots,a_5\leq 219$ is
a huge number: $(2\cdot219+1)^5=16305067506199\approx 1.63\cdot 10^{13}$.

We continue with the enumeration process.
In Case I we performed the enumeration
process with the following parameters:
\[
\begin{array}{|c|c|c|c|}\hline
                 & S  &  s  &  {\rm enumerated} \\ \hline
{\rm Step \;\;I} & 10^{258} & 10^{20} &  6 \\ \hline
{\rm Step \;\;II}& 10^{20}  & 10^{10}  & 6  \\ \hline
{\rm Step \;\;III}& 10^{10}   & 10^{3}  & 19922\\ \hline
{\rm Step \;\;IV}& 10^{3}   & 10^{2}  & 13506\\ \hline
{\rm Step \;\;V}& 10^{2}   &          & 1194\\ \hline
\end{array}
\]
All together the procedure took 2-3 minutes.
Parallely to enumeration we made sieving modulo 113, 787.
We found the only possible solution $a_1=\ldots=a_5=0$.

\noindent 
In Case 2 we have
\[
\frac{\frac{1}{2}X^{(i1)}}{\frac{1}{2}X^{(i2)}}=
\left(\frac{\varepsilon_1^{(i1)}}{\varepsilon_2^{(i1)}}\right)^{a_2-a_1}
\left(\frac{\varepsilon_3^{(i1)}}
{\varepsilon_1^{(i1)}\varepsilon_2^{(i1)}\varepsilon_4^{(i1)}}\right)^{a_4}
\left(\frac{\varepsilon_3^{(i1)}}
{\varepsilon_1^{(i1)}\varepsilon_2^{(i1)}\varepsilon_5^{(i1)}}\right)^{a_5}.
\]
Therefore in Case II we can determine the possible
values of $a_2-a_1,a_4,a_5$. We performed the enumeration algorithm
with the following parameters:
\[
\begin{array}{|c|c|c|c|c|}\hline
                 & S  & S^2 & s  &  {\rm enumerated} \\ \hline
{\rm Step \;\;I} & 10^{258} & 10^{516}  & 10^{20} &  0 \\ \hline
{\rm Step \;\;II}& 10^{20}  & 10^{40} & 10^{10}  & 0  \\ \hline
{\rm Step \;\;III}& 10^{10} & 10^{20}  & 10^{3}  & 38 \\ \hline
{\rm Step \;\;IV}& 10^{3}   & 10^6    &  10^{2}  & 202\\ \hline
{\rm Step \;\;V}& 10^{2}  & 10^4    &          & 79\\ \hline
\end{array}
\]
The procedure took all together some seconds.
In this case there was no way to diminish the possible exponent vectors by
sieving. For all the 319 possible values of $a_2-a_1,a_4,a_5$
we let $a_1,a_3$ run through the interval $[-219,219]$.
The exponent vectors $(a_1,\ldots,a_5)$ were tested modulo 
113, 787, 1223, 2053 if they satisfy the unit equation (\ref{u2}).

Finally we got three solutions of equation (\ref{u2}): 
\[
(a_1,\ldots,a_5)=(0,0,0,0,0),(0,1,0,0,0),(1,0,0,0,0).
\]
These correspond to 
\[
X=1,\;\; 1+\sqrt{\mu},\;\; 1-\sqrt{\mu},
\] 
that is
\[
\frac{1}{2}\cdot 1 +\frac{1}{2}\cdot 1 =1,\;\;
\frac{1}{2}\cdot (1+\sqrt{\mu}) +\frac{1}{2}\cdot (1-\sqrt{\mu}) =1,\;\;
\frac{1}{2}\cdot (1-\sqrt{\mu}) +\frac{1}{2}\cdot (1+\sqrt{\mu}) =1.
\]
By (\ref{ccc}) 
\[
X=\frac{\nu_G}{\nu_M},
\]
with $U=\nu_M$ and $U-2\sqrt{\mu}V=\nu_G$, we have
\[
\frac{V}{U}=\frac{1}{2\sqrt{\mu}}\left( 1-\frac{\nu_G}{\nu_M}   \right).
\]
This only gives and algebraic integer value for $X=1$ (out of the above
possible values of $X$), whence
the only solution of equation (\ref{Fex}) is $U=\nu_M,V=0$.

Following the general arguments of Section \ref{qqeq} we have
\[
Q_1(X,Y,Z)=X^2+\xi Z^2=U=\nu_M,\;\; Q_2(X,Y,Z)=Y^2-XZ=0.
\]
We get
\[
Q_0(X,Y,Z)=Y^2-XZ=0,
\]
with non-trivial solution $X_0=1,Y_0=0,Z_0=0$.
We set
\[
X=X_0R,\;\; Y=Y_0R+P,\;\; Z=Z_0R+Q,
\]
with parameters $P,Q,R\in M$. We obtain $P^2-RQ=0$. We multiply by $Q$ the
above equation and replace $RQ$ by $P^2$. Hence
\begin{eqnarray}
\kappa\cdot X&=&P^2,\nonumber\\
\kappa\cdot Y&=&PQ,\label{ss}\\
\kappa\cdot X&=&Q^2,\nonumber
\end{eqnarray}
and we replace $\kappa,P,Q$ by integer parameters $P,Q\in\Z_M$
(see Section \ref{qqeq}). It follows that $\kappa$ can only be a unit in $M$.
Substituting these representations into $Q_1(X,Y,Z)=U$
we obtain
\begin{equation}
P^4+\xi Q^4=\kappa^2\nu_M=\nu
\label{q1}
\end{equation}
where $\nu=\pm \varepsilon_1^{k_1},\varepsilon_2^{k_2}$ is a unit in $M$,
$\varepsilon_1,\varepsilon_2$ being the fundamental units in $M$.
Let $k_i=4k_i'+\ell_i$ with $-1\leq \ell_i\leq 2$ ($i=1,2$) and let
\[
P'=P\cdot \varepsilon_1^{-k_1'}\cdot \varepsilon_2^{-k_2'},
\;\;
Q'=Q\cdot \varepsilon_1^{-k_1'}\cdot \varepsilon_2^{-k_2'}.
\]
For $\nu'=\varepsilon_1^{\ell_1}\cdot \varepsilon_2^{\ell_2}$
we can easily solve equation
\[
(P')^4+\xi\cdot  (Q')^4=\pm \nu'
\]
in $P',Q'\in \Z_M$ using Theorem \ref{thth} and obtain that $P'=1,Q'=1$
is the only solution.
Therefore up to a unit factor in $M$ we have $(P,Q)=(1,0)$, 
whence up to a unit factor in $M$ we obtain
$(X,Y,Z)=(1,0,0)$. Hence up to equivalence the
only generator of relative power integral basis of $K$ over $M$ is $\xi$.

\section{Generators of absolute power integral bases}

As we have seen, up to equivalence 
$\xi$ is the only generator of relative power integral bases of $K$ over $M$.
By \cite{grsz} this implies that any generator of absolute power integral bases 
of $K$ must have the form
\[
\zeta=z_0+z_1\mu+z_2\mu^2 \pm \varepsilon_1^{k_1}\varepsilon_2^{k_2}\cdot \xi,
\]
where $z_0,z_1,z_2,k_1,k_2\in\Z$, $\varepsilon_1, \varepsilon_2$
are fundamental units in $M$.
We let $z_1,z_2,k_1,k_2$ run through the interval [-25,25].
In addition to $\xi$ we found only one algebraic integer of this shape
with index $<10^{15}$. This element has index 65329214857201.

\end{document}